\documentclass[11pt]{amsart}

\usepackage{amsmath,amssymb,amsthm,mathtools}
\usepackage{enumitem}
\usepackage{microtype}
\usepackage[hidelinks]{hyperref}

\newtheorem{theorem}{Theorem}[section]
\newtheorem{proposition}[theorem]{Proposition}
\newtheorem{lemma}[theorem]{Lemma}

\theoremstyle{definition}

\newtheorem{remark}[theorem]{Remark}

\newcommand{\Zz}{\mathcal{Z}}
\newcommand{\PP}{\mathbb{P}}
\newcommand{\ZZ}{\mathbb{Z}}
\newcommand{\PGL}{\operatorname{PGL}}
\newcommand{\Aut}{\operatorname{Aut}}
\newcommand{\Stab}{\operatorname{Stab}}

\title[Projective symmetries of Chebyshev polynomials]
      {Projective Symmetries of Chebyshev Polynomials}

\author{Saeed Tafazolian}
\author{Jaap Top}
\address{Departamento de Matem\'{a}tica - Instituto de Matem\'{a}tica, Estat\'{i}stica e Computação Cient\'{i}fica
(IMECC) - Universidade Estadual de Campinas (UNICAMP), Rua S\'{e}rgio Buarque de Holanda, 651, Cidade Universit\'{a}ria,  Zeferino Vaz, Campinas, SP 13083-859, Brazil}
\email{saeed@unicamp.br}
\address{Bernoulli Institute for Mathematics, Computer Science, and Artificial Intelligence\\
Nijenborgh~9\\9747 AG Groningen\\ the Netherlands}

\email{j.top@rug.nl}


\begin{document}

\begin{abstract}
Let $k$ be an algebraically closed field of odd characteristic $p$, and let
$\varphi_d$ be the monic Chebyshev polynomial defined by
$\varphi_d(t+t^{-1})=t^d+t^{-d}$. For $d>4$ with $p\nmid d$, let
$H_d\le\PGL_2(k)$ be the stabilizer of the zero set of $\varphi_d$. We show that
\[
H_d\cong
\begin{cases}
D_d,&\text{if }2d=p^r+1\;\text{ for some integer }\;r\geq 1,\\
S_3,&\text{if }4d=p^r+1\;\text{ for some integer }\; r\geq 1,\\
C_2,&\text{otherwise}.
\end{cases}
\]
\end{abstract}

\maketitle

\section{Introduction}\label{sec:Intro}
For an integer $d\ge1$, let $\varphi_d(X)\in\ZZ[X]$ be the monic Chebyshev
polynomial determined by
\begin{equation}\label{eq:cheb-def}
   \varphi_d(t+t^{-1})=t^d+t^{-d}.
\end{equation}
Thus $\varphi_0=2$, $\varphi_1=X$, and
$\varphi_{d+1}=X\varphi_d-\varphi_{d-1}$.  The polynomial
$\varphi_d$ is also known as the Dickson polynomial $D_d(X,1)$; see \cite{LMT} for
background on Dickson polynomials.

Let $k$ be an algebraically closed field of odd characteristic $p$, and assume
that $p\nmid d$. Denote by
\[
   \Zz_d=\{\alpha\in k:\varphi_d(\alpha)=0\}\subset\PP^1(k)
\]
the zero set of $\varphi_d$, and put
\[
   H_d=\Stab_{\Aut(\mathbb{P}^1)}(\Zz_d)=\left\{\tau\in \Aut(\mathbb{P}^1)\;:\;\tau(\Zz_d)=\Zz_d\right\}.
\]
Since $\varphi_d(-X)=(-1)^d\varphi_d(X)$, the involution
$\sigma:X\mapsto-X$ belongs to $H_d$.

For comparison in characteristic zero, Mosunov
\cite[Thm.~1.3]{Mosunov} computed the rational linear automorphism groups of
the homogenized Chebyshev forms of the first and second kinds.  For the first
kind, he obtains a cyclic group of order $2$ when the degree is odd and a
dihedral group $D_2$ of order $4$ when the degree is even; if transformations
sending the form to its negative are also allowed, the resulting group is
$D_2$ in both cases.  This is a different notion from the projective setwise
stabilizer considered here, but it gives a useful characteristic-zero
comparison.

 The main result proven in the present manuscript  is as follows.

\begin{theorem}\label{thm:stabilizer}
Let $k$ be algebraically closed of odd characteristic $p$ and let $d>4$ with
$p\nmid 2d$. Then
\[ \text{neither}\; 2d-1\; \text{nor}\; 4d-1 \;\text{is a power of}\;
p \Longleftrightarrow
H_d=\{X\mapsto X,\ X\mapsto-X\}.
\]
\end{theorem}
Denote by $D_d$ the dihedral group of order $2d$. The two exceptional cases
are described in Propositions~\ref{prop:Hd-2d} and~\ref{prop:Hd-4d}, namely
$H_d\cong D_d\;\text{if }2d=p^r+1$, and
   $H_d\cong S_3\;\text{if }4d=p^r+1$.
We consider the result and the techniques used in the proof interesting by
themselves; besides this, they can be used to describe the automorphism
groups of hyperelliptic curves defined by $y^2=\varphi_d(x)$ (at least when $d$
is even). In our manuscript \cite{TT} we use the knowledge of $H_d$ in every
relevant characteristic, to determine the automorphisms of the plane degree $d$
curves defined by $y^d=\varphi_d(x)$ for all pairs $(\text{char}(k), d)$ such that
$\text{char}(k)\nmid 2d$.

The paper is organized as follows. Section~2 recalls some basic facts about
$\varphi_d$, including a short elementary proof of its
separability. Section~3 contains two consequences of Lucas' theorem. In
Section~4 the stabilizer in $\Aut(\PP^1)$ of the set of zeros of
$\varphi_d$ is studied. The main results will be proven in
Section~5.

\section{Chebyshev polynomials}

\subsection{A derivative identity and separability}

A consequence of the relation \eqref{eq:cheb-def} defining $\varphi_d$ is the following.

\begin{lemma}\label{lem:derivative-identity}
For every $d\ge1$,
\begin{equation}\label{eq:derivative-identity}
   (X^2-4)\varphi_d'(X)^2
   =d^2\bigl(\varphi_d(X)^2-4\bigr)
\end{equation}
in $\ZZ[X]$.
\end{lemma}

\begin{proof}
The evaluation map
\[
   \ZZ[X]\longrightarrow\ZZ[t,t^{-1}],
   \qquad X\longmapsto t+t^{-1},
\]
is injective. Differentiating \eqref{eq:cheb-def} with respect to $t$ and
multiplying by $t$ gives
\[
   (t-t^{-1})\varphi_d'(t+t^{-1})
   =d(t^d-t^{-d}).
\]
Squaring this and applying the identities
\[
   (t-t^{-1})^2=(t+t^{-1})^2-4
\]
and
\[
   (t^d-t^{-d})^2=(t^d+t^{-d})^2-4
\]
completes the proof.
\end{proof}

Separability of $\varphi_d(X)\in k[X]$ whenever
$\text{char}(k)\nmid 2d$ is well known; see for
example \cite[Thm.~6.1(b)]{GS}, \cite[Thm.~7.2(a)]{KTT}, \cite[Lemma~2.5]{TTJNT}, \cite[Rem.~6.2]{ABS}. The preceding identity implies it as follows.

\begin{lemma}\label{lem:separable}
Assume that $k$ is algebraically closed of characteristic 
$p\nmid 2d$. Then $\varphi_d$ is separable over $k$. In particular the set $\Zz_d$ of zeros of $\varphi_d$ in $k$
satisfies
\[
   \#\Zz_d=d.
\]
\end{lemma}

\begin{proof}
Suppose that $\alpha$ were a multiple zero of $\varphi_d$. Then
$\varphi_d(\alpha)=\varphi_d'(\alpha)=0$. Considering
\eqref{eq:derivative-identity} over $k$, evaluation at $\alpha$ gives
\[
   0=-4d^2,
\]
which is impossible because $p$ is odd and $p\nmid 2d$. Hence $\varphi_d$ is
separable. Since $k$ is algebraically closed and $\deg\varphi_d=d$, the polynomial has
exactly $d$ distinct zeros in $k$.
\end{proof}

\subsection{Coefficients}
Write
\begin{equation}\label{eq:coeff}
   \varphi_d(X)=\sum_{j=0}^{\lfloor d/2\rfloor}c_jX^{d-2j},
   \qquad
   c_j=(-1)^j\frac{d}{d-j}\binom{d-j}{j}\in\ZZ.
\end{equation}
In particular,
\begin{equation}\label{eq:firstcoeff}
   c_0=1,\qquad c_1=-d.
\end{equation}

\subsection{The differential equation}
Differentiating \eqref{eq:derivative-identity} in $\ZZ[X]$ gives
\[
   2\varphi_d'\cdot
   \bigl((X^2-4)\varphi_d''+X\varphi_d'-d^2\varphi_d\bigr)=0.
\]
Since $\varphi_d'$ is not the zero polynomial, one obtains
\begin{equation}\label{eq:differential}
   (X^2-4)\varphi_d''+X\varphi_d'-d^2\varphi_d=0.
\end{equation}

\section{Lucas' theorem and a relation between certain polynomials}\label{sec:Lucas}
Recall the following result by \'{E}.~Lucas (1878).
\begin{theorem}[Lucas, \cite{Lucas}]\label{thm:lucas}
Let $p$ be prime and $N=\sum N_ip^i$ and $M=\sum M_ip^i$ with $0\le N_i,M_i<p$
the base $p$ expansions of $N,M\in\mathbb{Z}_{\geq 0}$.
Then $\binom NM\equiv\prod_i\binom{N_i}{M_i}\pmod p$; in particular
$\binom NM\not\equiv0\pmod p$ if and only if $M_i\le N_i$ for all $i$.
\end{theorem}

This will now be used to translate the property $2d\in 1+p^{\ZZ_{\geq 1}}$ into
a statement about a certain polynomial $P_d(Z)$ in characteristic $p$.
For every integer $d>0$ consider
\[
   P_d(Z):=\sum_{j=0}^{d}\binom {2d}{2j}Z^j\in k[Z],
\]
where $k$ is an algebraically closed field of characteristic $p\geq 3$.

\begin{lemma}[Scaling criterion]\label{lem:scaling}
Let $d\geq 2$ and $p\nmid 2d$, then
\[
P_d(aZ)=P_d(Z) \text{ for some }a\in k \text{ with } a\ne1 \Longleftrightarrow 2d\in 1+p^{\ZZ_{>0}}.\]
\end{lemma}

\begin{proof} $\Leftarrow$: Given is that $2d=1+p^r$ for some $r\geq 1$. Applying Theorem~\ref{thm:lucas}
with $N=2d$ and $M=2j\leq 2d$ yields $P_d(Z)=1+Z^{(1+p^r)/2}$. Choosing $a\in k$ of order dividing
$(1+p^r)/2$ and $a\neq 1$ completes the proof (note that the assumption that $k$ is algebraically
closed, or at least contains an algebraic closure of $\mathbb{F}_p$, is used here to guarantee that
$k$ contains the required roots of unity, for every $r\geq 1$).\\
$\Rightarrow$: 
Consider the base $p$ expansion $2d=\sum N_ip^i$ with $0\leq N_i<p$ for all $i$. We start by
showing that $N_i\in\{0,1\}$ for all $i$, and $N_0=1$.  Indeed, were $N_i\ge2$ for some $i\geq 0$ then
Theorem~\ref{thm:lucas} applied to $N=2d$ and $M=2p^i$ shows $\binom N{2p^i}\neq0$ in $k$, so $a^{p^i}=1$
which results in the contradiction
$a=1$. Thus every base $p$ digit of $2d$ is $0$ or $1$. Moreover, since $p\nmid 2d$ implies $N_0\neq 0$,
one concludes $N_0=1$.

Next, we show that besides $N_0$, exactly one more base $p$ digit of $2d$ equals $1$. As $2d$ is even and $p$ is odd, the number of base $p$ digits $N_r$ equal to $1$, is even. And were $N_r=N_s=1$ with $0<r<s$, then each of $1+p^r$, $1+p^s$ and $p^r+p^s$
are even and $<2d$. Theorem~\ref{thm:lucas} shows that the coefficients $\binom{2d}{1+p^r}, \binom{2d}{1+p^s}$
and $\binom{2d}{p^r+p^s}$ in $P_d\in k[Z]$
are nonzero, hence
$a^{(1+p^r)/2}=a^{(1+p^s)/2}=a^{(p^r+p^s)/2}=1$. 
 This implies
\[a^1=a^{(1+p^r)/2}\cdot a^{(1+p^s)/2}\cdot a^{-(p^r+p^s)/2}=1,\] a
contradiction. Hence  exactly one
higher base $p$ digit of $2d$ equals $1$ and $2d=1+p^r$.
\end{proof}

\begin{lemma}[Sparse criterion]\label{lem:sparse}
Let $p$ be prime and let $d\geq 2$ such that $p\nmid 2d$. If
$\binom {2d}{2j}\equiv0\pmod p$ for all $0<j<d$, then $2d=1+p^r$ for some $r\ge1$.
\end{lemma}

\begin{proof}
The hypothesis says that
\[
   P_d(Z)=1+Z^d.
\]
Since $d>1$ and $p\nmid d$, a nontrivial $d$-th root of unity
$a\in k^\times$ exists. Then $P_d(aZ)=P_d(Z)$, and Lemma~\ref{lem:scaling} shows
$2d=1+p^r$ for some $r\ge1$.
\end{proof}
\subsection*{Relating $\varphi_d$ and $P_d$}
The relevance of the polynomials $P_d$ towards our result, starts as follows.
\begin{lemma}\label{lem:coordinate}
\[
   (1-Z)^d\varphi_d\Bigl(2\frac{1+Z}{1-Z}\Bigr)=2P_{d}(Z).
\]
\end{lemma}

\begin{proof}
Introduce an indeterminate $R$ and substitute $Z=R^2$. Put
\[
T=\frac{1+R}{1-R},\;\text{ which implies }
T+T^{-1} = 2\frac{1+R^2}{1-R^2}.
\]
By \eqref{eq:cheb-def} it follows that
\[
\begin{aligned}
(1-Z)^d\varphi_d(2\frac{1+Z}{1-Z})=(1-R^2)^d
\varphi_d\left(2\frac{1+R^2}{1-R^2}\right)
&=(1-R^2)^d(T^d+T^{-d})\\
&=(1+R)^{2d}+(1-R)^{2d}\\
&=2\sum_{j=0}^{d}\binom{2d}{2j}R^{2j}\\
&=2P_{d}(R^2)=2P_d(Z),
\end{aligned}
\]
finishing the proof.
\end{proof}

\begin{proposition}\label{prop:fixpm2}
Let $k$ be algebraically closed of characteristic $p$ and $d\geq 5$
and $p\nmid 2d$. Then $H_d$ contains a nontrivial element fixing both $2$ and $-2$ if
and only if
$2d=1+p^r$
for some $r\ge1$.
\end{proposition}

\begin{proof}
Let
\[
   Z=\frac{X-2}{X+2}.
\]
Since $\varphi_d(2)=2$ and $\varphi_d(-2)=2(-1)^d$, neither $2$ nor $-2$ belongs to $\Zz_d$.
The automorphism $X\mapsto Z$ on $\PP^1$ sends $2$ to $0$ and $-2$ to $\infty$. By
Lemma~\ref{lem:coordinate}, it maps $\Zz_d$ onto the zero set of
$P_{d}$.

Suppose first that $H_d$ contains a nontrivial element fixing both $2$ and
$-2$. In the $Z$-coordinate it then fixes $0$ and $\infty$, hence is given by $Z\mapsto aZ$ with
$a\in k^\times$, $a\ne1$. Hence $P_{d}(aZ)=cP_{d}(Z)$ for some
$c\in k^\times$. Since $P_{d}(0)=1$, we have $c=1$, and
Lemma~\ref{lem:scaling} gives $2d=1+p^r$.

Conversely, assume $2d=q+1$ with $q=p^r$. Introducing $R$ with $Z=R^2$ one obtains
\[
\begin{aligned}
 2P_d(Z)=(1+R)^{q+1}+(1-R)^{q+1}
 &= (1+R)(1+R^q)+(1-R)(1-R^q)\\
 &=2(1+R^{q+1})=2(1+Z^d),
\end{aligned}
\]
hence $P_{d}(Z)=1+Z^d$.
Now
$Z\mapsto\zeta Z$, for a nontrivial $d$-th root of 
unity $\zeta$, preserves the zero set of $P_d$, hence 
defines an element $\neq\textit{id}$ of $H_d$ fixing $2$ and $-2$.
\end{proof}

\section{Finite projective stabilizers}
In this section $k$ denotes an algebraically closed
field of characteristic $p\geq 0$, and $d\geq 3$ is
an integer. Throughout, it is assumed that
$p\nmid 2d$.
By Lemma~\ref{lem:separable}, these conditions imply
that the set $\Zz_d$ of zeros of the Chebyshev polynomial
$\varphi_d$ in $k$ satisfies $\#\Zz_d=d\ge3$. We study the
group
\[H_d:=\{\tau\in \Aut(\PP^1)\;:\;\tau(\Zz_d)=\Zz_d\}.\]
Restriction to
$\Zz_d$ gives a homomorphism
\[
   H_d\longrightarrow\operatorname{Sym}(\Zz_d).
\]
Elements in its kernel fix every point of $\Zz_d$, hence at least three distinct points
of $\PP^1(k)$. Thus the homomorphism is injective and $H_d$ is a finite group. As already observed in 
Section~\ref{sec:Intro}, the equality $\varphi_d(-X)=(-1)^d\varphi_d(X)$ shows that
\[\sigma\colon X\mapsto -X\]
is an element of order $2$ in $H_d$. 
\begin{remark}
If $d=3$ then $\varphi_d(X)=X^3-3X$. In every characteristic $\neq 3$, all permutations of the
three zeros of $\varphi_3$ lift to elements of $H_3$. Hence here
\[H_3\cong S_3.\]
For $d=4$ one finds $\varphi_4(X)=X^4-4X^2+2$.
A suitable scaling $Y=\lambda X$ transforms $\frac12\varphi_4$ into $Y^4+aY^2+1$ with $a^2=8$.
A calculation (conveniently done using cross ratios, see for example \cite[Section~6.5.2]{Dolgachev})
results in
\[ H_4\cong\left\{
\begin{array}{rl}
\ZZ/2\ZZ\times \ZZ/2\ZZ & \text{if char}(k)\neq 5,7;\\
A_4 &  \text{if char}(k)=5;\\
D_4 &   \text{if char}(k)=7.
\end{array}
\right.\]
\end{remark}

\begin{lemma}\label{lem:second-involution}
If $H_d$ has no involution other than $\sigma$, then
$H_d=\langle\sigma\rangle$.
\end{lemma}

\begin{proof}
Assume that $\sigma$ is the unique involution of $H_d$. For every $g\in H_d$,
the conjugate $g\sigma g^{-1}$ is again an involution, hence it equals
$\sigma$. Thus
\[
   H_d\le C_{\Aut(\PP^1)}(\sigma),
\]
the centralizer in $\Aut(\PP^1)$ of $\sigma$. This
centralizer consists of the maps
\[
   X\mapsto aX
   \qquad\text{and}\qquad
   X\mapsto a/X,
   \qquad a\in k^\times.
\]
Every map of the second type is an involution, so by
assumption is not in
$H_d$. Hence every $g\in H_d$ has the form
$X\mapsto aX$.

Such $g$ preserves $\Zz_d$, hence the monic and separable polynomials $\varphi_d(X)$ and
$a^{-d}\varphi_d(aX)$ are equal, hence
\[
   \varphi_d(aX)=a^d\varphi_d(X).
\]
Comparing the coefficients of $X^{d-2}$ and using \eqref{eq:firstcoeff} gives
$a^2=1$, because $p\nmid d$. Therefore  $H_d=\langle\sigma\rangle$.
\end{proof}

Assume now that
$H_d\supsetneq\langle\sigma\rangle$. By Lemma~\ref{lem:second-involution},
$H_d$ contains an involution $\tau\ne\sigma$. Put
\[
   \mu:=\sigma\tau\in H_d.
\]
Then $\mu\in H_d$ and
\begin{equation}\label{eq:reversible}
   \sigma \mu\sigma=\mu^{-1}.
\end{equation}
Note that $\mu\neq \textit{id}$ since $\sigma\neq \tau$,
and $\mu\neq \sigma$ since $\tau\neq \textit{id}$.

\begin{lemma}[Normal form]\label{lem:normalform}
There exist $u\in k$ and $s\in k^\times$ such that
\begin{equation}\label{eq:normalform}
   \mu(X)=u+\frac{s}{X+u},
\end{equation}
and $\mu$ is an involution if and only if $u=0$.
\end{lemma}

\begin{proof}
Write $u:=\mu(\infty)\in\PP^1$.
If $u=\infty$, then $\mu(X)=aX+b$ for certain $a,b\in k$. Since $\tau(X)=-aX-b$ is an involution, we have
$a^2=1$ and $(a-1)b=0$.
If $a=-1$, then $b=0$ because the characteristic
is not $2$, and hence $\tau=\mathrm{id}$,
a contradiction. Thus $a=1$. Hence
$\mu(X)=X+b$.
Since $\mu\in H_d$, the polynomials $\varphi_d(X+b)$ and $\varphi_d(X)$
are monic of degree $d$ and separable, with the same zero set. Hence the two polynomials
are equal. Comparing the coefficients of $X^{d-1}$
gives $db=0$, implying $b=0$, a contradiction. Thus $u\neq \infty$, i.e., $u\in k$.

The equality $\sigma \mu\sigma=\mu^{-1}$ and the fact that
$\mu^{-1}(u)=\infty$ show that $\mu(-u)=\infty$.
Hence $\mu$ is given as
\[ \mu(X) = u + \frac{s}{X+u}\]
for some $s\in k^\times$.
From $\mu(\mu(-u))=u$ it follows that $\mu$ is an
involution if and only if $u=0$.
\end{proof}

\subsection{The root identity}

Let $\mu\in H_d$ be as before, given as in \eqref{eq:normalform}.
Since $\mu^{-1}$ preserves $\Zz_d$, we have $\varphi_d(u)\ne0$: indeed, the proof of Lemma~\ref{lem:normalform}
shows $\mu^{-1}(u)=\infty\not\in\Zz_d$.

Write the $d$ zeros of $\varphi_d$ as $\alpha_1,\dots,\alpha_d$. Since $\mu^{-1}$ permutes them
and $\mu^{-1}(\alpha_i)=-u-\frac{s}{u-\alpha_i}$, one computes
\[
\begin{aligned}
(X+u)^d\varphi_d\bigl(\mu(X)\bigr)
&=\prod_{i=1}^d\bigl[(X+u)(u+\frac{s}{X+u}-\alpha_i)\bigr]\\
&=\prod_{i=1}^d\bigl[(u-\alpha_i)(X+u)+s\bigr]\\
&=\varphi_d(u)\prod_{i=1}^d\bigl(X-\mu^{-1}(\alpha_i)\bigr)\\
&=\varphi_d(u)\varphi_d(X),
\end{aligned}
\]
i.e.,
\begin{equation}\label{eq:rootidentity}
   (X+u)^d \varphi_d\Bigl(u+\frac{s}{X+u}\Bigr)=\varphi_d(u)\varphi_d(X).
\end{equation}

\subsection{A logarithmic-differential identity}

Set
\begin{equation}\label{eq:Adef}
   A(X):=-\frac{X}{X^2-4},
   \qquad
   \omega:=A(X)\,\text{d}X=-\frac{X\,\text{d}X}{X^2-4}.
\end{equation}
Thus $\omega=-\tfrac12 \text{d}\log(X^2-4)$.  The differential $\omega$ has simple
poles at $2$ and $-2$, both of residue $-1/2$, and a simple pole at $\infty$
of residue $1$.

\begin{lemma}[Logarithmic-differential identity]\label{lem:root-connection}
Assume that $k$ has characteristic $p\geq 0$, that $d\geq 5$, and that $p\nmid 2d$.
Let
\[
   \mu(X)=u+\frac{s}{X+u}\in H_d,
   \qquad s\ne0.
\]
Then
\begin{equation}\label{eq:connection}
   A(X)=\frac{2(d-1)}{X+u}+\mu'(X)A(\mu(X))
\end{equation}
as an identity in $k(X)$.  Equivalently (multiplying by $\text{\em d}X$), one has
\begin{equation}\label{eq:connection-differential}
   \omega
   =2(d-1)\frac{\textrm{\em d}X}{X+u}+\mu^*\omega.
\end{equation}
\end{lemma}

\begin{proof}
Using \eqref{eq:rootidentity}, put
\[
   H_\mu(X):=(X+u)^d \varphi_d(\mu(X))=\varphi_d(u)\varphi_d(X).
\]
Let $\beta$ be a zero of $\varphi_d$.  By Lemma~\ref{lem:separable},
$\varphi_d'(\beta)\ne0$.  Moreover $\beta\ne\pm2$, since
$   \varphi_d(2)=2$ and $\varphi_d(-2)=2(-1)^d$.
Also, $\beta\ne-u$ because $\mu$ sends $-u$ to $\infty$ whereas
$H_d$ preserves $\Zz_d\ni \beta$. The same remarks apply with $\beta$ replaced by $\mu(\beta)$.

Evaluating the differential equation \eqref{eq:differential} at a zero $\alpha$ of $\varphi_d$ gives
\begin{equation}\label{eq:root-logder}
   \frac{{\varphi_d}''(\alpha)}{{\varphi_d}'(\alpha)}
   =-\frac{\alpha}{\alpha^2-4}
   =A(\alpha).
\end{equation}
At $X=\beta$, differentiation of $H_\mu(X)=(X+u)^d\varphi_d(\mu(X))$ gives
\[
   \frac{H_\mu''(\beta)}{H_\mu'(\beta)}
   =\frac{2d}{\beta+u}
    +\mu'(\beta)\frac{{\varphi_d}''(\mu(\beta))}{{\varphi_d}'(\mu(\beta))}
    +\frac{{\mu}''(\beta)}{\mu'(\beta)}.
\]
As $H_\mu=\varphi_d(u)\varphi_d$ and
   $\mu'(X)=-{s}/{(X+u)^2}$ and 
   ${\mu''(X)}/{\mu'(X)}=-2/{(X+u)}$,
from \eqref{eq:root-logder} one obtains
\begin{equation}\label{eq:connection-on-roots}
   A(\beta)=\frac{2(d-1)}{\beta+u}
             +\mu'(\beta)A(\mu(\beta)).
\end{equation}
Put $t=X+u$, $B=ut+s$, and $Q=B^2-4t^2$.
Since
\[
\mu(X)=\frac{B}{t},
\qquad
\mu'(X)=-\frac{s}{t^2},
\]
we have
\[
\mu'(X)A(\mu(X))=\frac{sB}{tQ}.
\]
Consequently,
\[
\begin{aligned}
&A(X)-\frac{2(d-1)}{X+u}-\mu'(X)A(\mu(X))\\
&\qquad=
\frac{-XtQ-(X^2-4)\bigl(2(d-1)Q+sB\bigr)}
     {(X^2-4)tQ}.
\end{aligned}
\]
The numerator has degree at most $4$.
For every $\beta\in \Zz_d$, the denominator is nonzero:
indeed, $\beta\ne\pm2$, $\beta\ne-u$, and
$\mu(\beta)\in \Zz_d$ implies $\mu(\beta)\ne\pm2$,
so $Q(\beta)\ne0$.
By \eqref{eq:connection-on-roots}, the numerator therefore vanishes at all
$d\ge5$ distinct points of $\Zz_d$.
Hence it is identically zero, proving \eqref{eq:connection}.

\end{proof}

\begin{remark}\label{rem:connection-residue}
The differential form \eqref{eq:connection-differential} makes the local
information transparent. Since $\mu(-u)=\infty$ and $\mu$ is a M\"obius transformation, residues are
preserved under pullback. Indeed, if $t$ is a local parameter at a point $P$
and $z=c t+O(t^2)$, $c\ne0$, is a local parameter at $\mu(P)$, then the
coefficient of $dt/t$ is unchanged. Hence
\[
   \operatorname{res}_{-u}(\mu^*\omega)
   =\operatorname{res}_{\infty}(\omega)=1.
\]
Thus the total residue of the right-hand side at $X=-u$ is $2d-1$.
Furthermore, whenever $\mu^{-1}(2)$ or $\mu^{-1}(-2)$ is finite, it is a genuine
simple pole of $\mu^*\omega$, with residue $-1/2$.   Neither of these points equals $-u$, since $\mu(-u)=\infty$.
Hence these poles cannot be cancelled by the term
$2(d-1)\,dX/(X+u)$, whose only possible finite pole
is at $X=-u$.
\end{remark}

\begin{lemma}[The $4d$ Frobenius criterion]\label{lem:four-d-criterion}
Assume that $k$ is algebraically closed of  characteristic $p$, that $d\geq 1$,
and that $p\nmid 2d$. Put
\[
   h(X):=\frac{2X+12}{2-X},
   \qquad
   F_d(X):=\frac{(2-X)^d}{2}\,\varphi_d(h(X))\in k[X].
\]
Then $F_d=\varphi_d$ if and only if
$4d=1+p^r$
for some $r\ge1$.
\end{lemma}

\begin{proof}
Introduce $Y$ with $Y^2=X+2$ and put
\[
   T=\frac{(Y+2)^2}{2-X}=\frac{Y+2}{2-Y}.
\]
Then $T+T^{-1}=h(X)$, and the defining identity for $\varphi_d$ gives
\[
   2F_d(X)=(Y+2)^{2d}+(Y-2)^{2d}.
\]
Now substitute
\[
   X=Z^4+Z^{-4},
   \qquad
   Y=Z^2+Z^{-2}.
\]
Since $Y+2=(Z+Z^{-1})^2$ and $Y-2=(Z-Z^{-1})^2$, we obtain
\[
   2F_d(Z^4+Z^{-4})
   =(Z+Z^{-1})^{4d}+(Z-Z^{-1})^{4d}.
\]
On the other hand,
\[
   2\varphi_d(Z^4+Z^{-4})=2(Z^{4d}+Z^{-4d}).
\]
The substitution $k[X]\hookrightarrow k[Z,Z^{-1}]$, $X\mapsto Z^4+Z^{-4}$,
is injective. Therefore $F_d=\varphi_d$ if and only if
\[
   (Z+Z^{-1})^{4d}+(Z-Z^{-1})^{4d}=2(Z^{4d}+Z^{-4d}).
\]
Expanding the left-hand-side shows that this is equivalent to
\[
   \binom {4d}{2i}\equiv0\pmod p
   \qquad(0<i<2d).
\]
Applying Lemma 3.3 with $2d$ in place of $d$,
we obtain $4d=1+p^r$ for some $r\ge1$.

Conversely, if $4d=1+q$ with $q=p^r$, then one computes
\[
\begin{aligned}
 (Z+Z^{-1})^{q+1}+(Z-Z^{-1})^{q+1}
 &= (Z+Z^{-1})(Z^q+Z^{-q})
 +(Z-Z^{-1})(Z^q-Z^{-q})\\
 &=2(Z^{q+1}+Z^{-q-1}),
\end{aligned}
\]
showing $F_d=\varphi_d$.
\end{proof}

\subsubsection{The case $u=0$}

\begin{proposition}\label{prop:u0}
Assume that $k$ has characteristic $p$, that $d\geq 5$ and $p\nmid 2d$. \\ If $u=0$,
then
   $2d=1$ {in } $k$ and
   $s=\pm4$
and $H_d$ contains a nontrivial element fixing both $2$ and $-2$.
\end{proposition}

\begin{proof}
When $u=0$, the differential $\omega$ is regular at $X=0$. By
Remark~\ref{rem:connection-residue}, the residue of the right-hand-side of
\eqref{eq:connection-differential} at $X=0$ is $2d-1$. Hence $2d=1$ {in }$k$.
Both terms on the right in \eqref{eq:connection-differential} have at most simple poles at $0$, so this residue
cancellation makes the singularity at $0$ removable. Note
that $u=0$ means that $\mu$ is given by $\mu(X)=s/X$.
The remaining finite poles of $\mu^*\omega$ are the genuine simple poles
\[
   \mu^{-1}(2)=s/2,
   \qquad
   \mu^{-1}(-2)=-s/2,
\]
each of residue $-1/2$. The left-hand-side $-X\text{d}X/(X^2-4)$ has finite poles precisely at
$2$ and $-2$. Therefore
   $\{s/2,-s/2\}=\{2,-2\}$,
so $s=\pm4$. Hence one of $\mu$ and $\sigma \mu$ equals
$X\longmapsto\frac4X$,
which is nontrivial and fixes $2$ and $-2$.
\end{proof}

\subsubsection{The endpoint case $u^2=4$}

Since $\mu^{-1}(X)=-u+s/(X-u)$,  replacing $\mu$ by $\mu^{-1}$ changes $u$ to $-u$ and leaves $s$ unchanged.
Thus we may assume $u=2$.

\begin{proposition}\label{prop:endpoint}
Assume that $k$ has characteristic $p$ and $d\geq 5$ and $p\nmid 2d$, and suppose
$u^2=4$. Then, after replacing $\mu$ by $\mu^{-1}$ if necessary,
$4d=1, s=-16$ {in }$k$,
so $\mu^{-1}(X)=(2X+12)/(2-X)$. 

In $\ZZ$, it follows that
$4d=1+p^r$
for some $r\ge1$.
\end{proposition}

\begin{proof}
Take $u=2$. The differential $\omega$ has residue $-1/2$ at $X=-2$.
By Remark~\ref{rem:connection-residue}, the residue of the right-hand side of
\eqref{eq:connection-differential} at $-2$ is $2d-1$. 
In $k$ this implies
$2d-1=-\frac12$,
hence $4d=1$.

Since $\mu(\infty)=u=2$, the inverse image of the pole $2$ of $\omega$ is
$\infty$ and $\mu^{-1}(-2)$ is finite and is a genuine simple
pole of $\mu^*\omega$ of residue $-1/2$. It cannot equal $-2$ because $\mu(-2)=\infty$. Thus comparison with the two finite poles $\{2,-2\}$ of $\omega$
forces
$\mu^{-1}(-2)=2$,
or equivalently $2+s/4=\mu(2)=-2$, implying $s=-16$.
Therefore
\begin{equation}\label{eq:M3}
   \mu^{-1}(X)=-2+\frac{-16}{X-2}
            =\frac{2X+12}{2-X}.
\end{equation}
Applying \eqref{eq:rootidentity} to $\mu^{-1}$ gives
\[
   (X-2)^d \varphi_d(\frac{2X+12}{2-X})=\varphi_d(-2)\varphi_d(X)=2(-1)^d \varphi_d(X).
\]
Since $(2-X)^d=(-1)^d(X-2)^d$, this is equivalent to
\[
   \frac{(2-X)^d}{2}\varphi_d(\frac{2X+12}{2-X})=\varphi_d(X).
\]
Lemma~\ref{lem:four-d-criterion} therefore yields
$4d=1+p^r$
for some $r\ge1$.
\end{proof}

\subsubsection{The generic case}

\begin{proposition}\label{prop:generic}
Assume that $k$ has characteristic $p$ and $d\geq 5$ and $p\nmid 2d$ and $u\not\in\{0,\pm 2\}$.
Then
$2d=1,\;s=4-u^2$ in $k$.
In particular, $\mu$ fixes both $2$ and $-2$.
\end{proposition}

\begin{proof}
Since $u^2\ne4$, the differential $\omega$ is regular at $X=-u$. By
Remark~\ref{rem:connection-residue}, the residue of the right-hand side of
\eqref{eq:connection-differential} is $2d-1$. Hence
$2d=1$ in $k$.

The two terms on the right in \eqref{eq:connection-differential} have at most simple poles at $-u$, so the
singularity there is removable. Moreover $\mu(\infty)=u\notin\{2,-2\}$,
so both $\mu^{-1}(2)$ and $\mu^{-1}(-2)$ are finite. By
Remark~\ref{rem:connection-residue} they are distinct genuine simple poles of
$\mu^*\omega$, each with residue $-1/2$. Comparison with the finite
poles of $\omega$ gives
\[
   \mu^{-1}\{2,-2\}=\{2,-2\}.
\]
If $\mu(2)=-2$ then
   $\mu(-2)=2$
and $(u+2)^2=-s=(u-2)^2$.
Thus $8u=0$, impossible because the characteristic is odd and $u\ne0$.
Hence $\mu$ fixes both $2$ and $-2$, implying 
$s=4-u^2$.
\end{proof}

\section{Proof of the main theorem}\label{sec:main-proofs}

\begin{proof}[Proof of Theorem~\ref{thm:stabilizer}]
Recall that $\sigma:X\mapsto -X$ belongs to $H_d$,
and assume that
$H_d\supsetneq\langle\sigma\rangle$.
By Lemma~\ref{lem:second-involution}, $H_d$ contains an involution
$\tau\ne\sigma$; put
$\mu:=\sigma\tau$ so that
$\mu\ne\textit{id}$ and
   $\sigma \mu\sigma=\mu^{-1}$.
By Lemma~\ref{lem:normalform}, there exist $u\in k$ and $s\in k^\times$ such
that
$\mu(X)=u+{s}/{(X+u)}$.
We distinguish three cases.
\begin{itemize}
\item If $u=0$, then Proposition~\ref{prop:u0} shows that $H_d$ contains a
nonidentity element fixing both $2$ and $-2$. Hence
Proposition~\ref{prop:fixpm2} gives
$2d=1+p^r$
for some $r\ge1$.
\item
If $u^2=4$, Proposition~\ref{prop:endpoint} shows
$4d=1+p^r$
for some $r\ge1$.
\item
Finally, suppose that $u\ne0$ and $u^2\ne4$. Then
Proposition~\ref{prop:generic} shows
$2d=1, s=4-u^2$ {in } $k$,
and $\mu\neq\textit{id}$ fixes both $2$ and $-2$. Thus Proposition~\ref{prop:fixpm2} gives
$2d=1+p^r$
for some $r\ge1$.
\end{itemize}
Therefore
\[
H_d\supsetneq\langle\sigma\rangle
\quad\Longrightarrow\quad
2d=p^r+1
\ \text{or}\
4d=p^r+1.
\]

We now prove the converse. If $2d=p^r+1$, Proposition~\ref{prop:fixpm2}
produces a nonidentity element of $H_d$ fixing $2$ and $-2$. Since such an element
is not in $\langle\sigma\rangle$, one concludes $H_d\supsetneq\langle\sigma\rangle$.

If $4d=p^r+1$, Lemma~\ref{lem:four-d-criterion} implies that $\mu$ defined by
   $\mu(X):=(2X+12)/(2-X)$
is in $H_d$. Hence also in this case $H_d\supsetneq\langle\sigma\rangle$.

This completes the proof.
\end{proof}

\begin{remark}\label{rem:charzero}
The same argument also gives the characteristic-zero case. If $k$ is
algebraically closed of characteristic zero and $d>4$, then
$H_d=\langle\sigma\rangle$.
Indeed, a second involution would lead, according as $u=0$, $u^2=4$, or
$u\ne0$ and $u^2\ne4$, to the equality $2d=1$ or $4d=1$ in $k$, which is
impossible in characteristic zero.
\end{remark}

Denote by $D_d$ the dihedral group of order $2d$.
\begin{proposition}\label{prop:Hd-2d}
Assume that $k$ has characteristic $p$ and $d\geq 5$ and $2d=p^r+1$ for some $r\ge1$. Then
$H_d\cong D_d$.
\end{proposition}

\begin{proof}
Put $q=p^r$. By Lemma~\ref{lem:coordinate} and the proof of
Proposition~\ref{prop:fixpm2}, the coordinate
$Z=(X-2)/(X+2)$
identifies $\Zz_d$ with the zero set of $1+Z^d$.  Choose $\eta\in k^\times$
with $\eta^d=-1$.  After the further scaling $Z\mapsto\eta^{-1}Z$, one obtains that $H_d$
 is conjugate to
\[
   K:=\Stab_{\Aut(\PP^1)}(\mu_d),
   \qquad
   \mu_d=\{z\in k^\times:z^d=1\}.
\]
The group $K$ contains
\[
   R:=\{Z\mapsto\zeta Z\colon \zeta^d=1\}\cong C_d
\]
and the involution $\iota:Z\mapsto Z^{-1}$. Together
they generate a subgroup
$\cong D_d$ of $K$.
We now show that
$K$ is exactly this dihedral subgroup. 

First, $p\nmid \#K$, since
otherwise an element $g\in K$ of order $p$ exists.  Such an
element has exactly one fixed point in $\PP^1(k)$; consequently
its action on the $d$-element set $\mu_d$ has orbits of length $p$ in addition to at
most one fixed point. Thus
$d\equiv0\;\text{or}\;1\pmod p$.
On the other hand $2d=q+1$ gives
$2d\equiv1\pmod p$,
 contradicting the above.
Hence $p\nmid \#K$.

By the classification (originally due to L.E.~Dickson) of finite subgroups of $\PGL_2$
\cite[Thm.~A.8]{HKT}, the noncyclic group $K$ is either dihedral or one of
$A_4,S_4,A_5$.  Since $K$ contains an element of order $d\ge5$, the groups
$A_4$ and $S_4$ are impossible.  If $K\cong A_5$, the existence of an element of
order $d\ge5$ in $K$ forces $d=5$. But then $q=2d-1=9$, so $p=3$, contradicting
$p\nmid|A_5|$.
Hence $K$ is dihedral, say of order $2n$.

Let $C_n$ be its cyclic subgroup of index two. Since every element outside
$C_n$ has order $2$ and $R$ contains an element of order $d\ge5$, we have
$R\le C_n$.  Thus $C_n$ centralizes $R$.  A generator $Z\mapsto\zeta Z$ of
$R$ has $\zeta\ne\pm1$, and its centralizer in $\PGL_2(k)$ consists of the
scalings $Z\mapsto aZ$.  Therefore every element of $C_n$ is a scaling.
A scaling preserving $\mu_d$ must satisfy $a^d=1$, so $C_n=R$.  Hence
$n=d$ and $K\cong D_d$, finishing the proof.
\end{proof}

\begin{proposition}\label{prop:Hd-4d}
Assume that $k$ has characteristic $p$ and $d>4$ and $4d=p^r+1$ for some $r\ge1$. Then
$H_d\cong S_3$.
\end{proposition}

\begin{proof}
The proof of Theorem~\ref{thm:stabilizer} constructs
$\mu(X)=(2X+12)/(2-X)$
and shows that the given relation between $d$ and $p$
implies $\mu\in H_d$. A straightforward calculation reveals that
\[
   \langle\sigma,\mu\rangle\cong S_3\le H_d.
\]
It remains to prove that there are no further elements
in $H_d$.
We first determine all involutions in $H_d$. Let $\tau\ne\sigma$ be an
involution and put $\nu=\sigma\tau$; write
$\nu(X)=u+\frac{s}{X+u}$
as in Lemma~\ref{lem:normalform}. Since $4d=1$ in $k$, the cases $u=0$ and
$u\not\in\{0,\pm2\}$ are impossible by Propositions~\ref{prop:u0} and
\ref{prop:generic}. Thus
$u^2=4$. Proposition~\ref{prop:endpoint} then shows, after replacing $\nu$ by
$\nu^{-1}$ if necessary, that $\nu^{-1}=\mu$. 
Consequently the involutions in $H_d$ are exactly
$\{\sigma, \sigma \nu=\tau, \sigma \nu^{-1}=\sigma\tau\sigma\}$.

Conjugation on this three-element set defines a homomorphism
\[
   \rho\colon H_d\longrightarrow S_3.
\]
We claim that $\rho$ is injective. Indeed, let $g\in\ker\rho$. Then $g$ commutes with
$\sigma$, so $g$ has the form $X\mapsto aX$ or $X\mapsto a/X$. In the second
case $g$ is an involution. Since the three involutions above are all the
involutions in $H_d$ and neither $\tau$ nor $\sigma\tau\sigma$ commutes with
$\sigma$, this is impossible. Hence $g(X)=aX$. As in
the proof of Lemma~\ref{lem:second-involution},  comparison of
the coefficient of $X^{d-2}$ in $\varphi_d(aX)$ and
$a^d\varphi_d(X)$ shows $a^2=1$. Thus $g$ is either $\textit{id}$ or
$\sigma$. Since clearly $\sigma\notin\ker\rho$,
this shows that $\rho$ is injective.
Since $H_d$ contains the subgroup
$\langle\sigma,\mu\rangle\cong S_3$, it follows that $H_d\cong S_3$.
\end{proof}

\section*{Acknowledgments}
The first author was partially supported by the Brazilian funding agencies
CNPq (grant 302774/2025-4), FAPESP (grant 2024/00923-6) and FAEPEX
(grant 3485/25).

\section*{Declaration on the use of generative AI}
 The problem addressed in this paper arose from the authors’ earlier work on automorphisms of Chebyshev plane curves. The first author used Claude Opus 5.5 as an aid in developing the argument presented in Lemma 4.4. The authors subsequently reviewed and revised the argument and take full responsibility for the mathematical content of the paper.

\end{document}